\documentclass[a4paper,11pt]{article}
\usepackage{hyperref} 
\usepackage{amsfonts,latexsym,amstext}
\usepackage{amsmath}
\usepackage{amssymb}
\usepackage[english]{babel}
\usepackage[latin1]{inputenc}
\usepackage{amsthm}
\usepackage{latexsym}
\usepackage{array}
\usepackage{mathrsfs}
\usepackage{amsxtra}
\usepackage{amscd}
\usepackage{mathtools}
\usepackage{verbatim}
\usepackage{enumerate}

\usepackage[usenames]{color}
\definecolor{beige}{rgb}{0.96, 0.96, 0.86}
\definecolor{airforceblue}{rgb}{0.36, 0.54, 0.66}
\definecolor{antiquefuchsia}{rgb}{0.57, 0.36, 0.51}
\definecolor{awesome}{rgb}{1.0, 0.13, 0.32}

\usepackage{xcolor}

\usepackage{fancyhdr}

\newcommand{\tr}{\mathop{\mathrm{Tr}}\nolimits}

\renewcommand{\tilde}{\widetilde}

\definecolor{red}{rgb}{1.0,0.0,0.0}

\definecolor{blu}{rgb}{0.0,0.0,1.0}

\definecolor{gre}{rgb}{0.03,0.50,0.03}

\definecolor{darkviolet}{rgb}{0.58, 0.0, 0.83}

\usepackage{enumerate}

\newtheorem{theorem}{Theorem}[section]

\newtheorem{proposition}[theorem]{Proposition}

\newtheorem{assumption}[theorem]{Assumption}

\newtheorem{corollary}[theorem]{Corollary}

\numberwithin{equation}{section}

\def\qed{{\hfill\hbox{\enspace${ \square}$}} \smallskip}
\def\sqr#1#2{{\vcenter{\vbox{\hrule height .#2pt \hbox{\vrule
 width .#2pt height#1pt \kern#1pt \vrule
width .#2pt} \hrule height .#2pt}}}}
\def\square{\mathchoice\sqr54\sqr54\sqr{4.1}3\sqr{3.5}3}

\def\epsilon{\varepsilon}

\def\ds{\begin{displaystyle}}
\def\eds{\end{displaystyle}}

\def\<{\left\langle }
\def\>{\right\rangle }

\def\dim{\noindent \hbox{{\bf Proof.} }}

\def\R{\mathbb R}

\def\E{\mathbb E}

\def\cald{{\cal D}}

\title{An optimal advertising model with carryover effect
\\
and mean field terms}
\author{Fausto Gozzi
\\
Dipartimento di Economia e Finanza\\
Universit\`a LUISS - Guido Carli\\
Viale Romania 32,
00197 Roma,
Italy\\
e-mail: fgozzi@luiss.it\\
\\
Federica Masiero\\
Dipartimento di Matematica e Applicazioni\\ Universit\`a di Milano Bicocca\\
via Cozzi 55, 20125 Milano, Italy\\
e-mail: federica.masiero@unimib.it\\
\\
Mauro Rosestolato\\
Dipartimento di Economia\\ Universit\`a di Genova \\
via Vivaldi 5, 12126 Genova, Italy\\
e-mail: mauro.rosestolato@unige.it}

\date{}

\makeatletter
\newcases{mycases}{\quad}{\hfil$\m@th\displaystyle{##}$}{$\m@th\displaystyle{##}$\hfil}{\lbrace}{.}
\makeatother

\DeclareMathOperator*{\argmax}{arg\,max}

\usepackage{cancel}

\begin{document}
\maketitle

\begin{abstract}
  We consider a class of optimal advertising problems under uncertainty for the introduction of a new product into the market, on the line of the seminal papers of Vidale and Wolfe, 1957 (\cite{VW}) and Nerlove and Arrow, 1962 (\cite{NA}).
  The main features of our model are that, on one side, we assume a carryover effect
  (i.e.\ the advertisement spending affects the
  goodwill with some delay); on the other side
  we introduce, in the state equation and in the objective, some mean field terms that take into account the presence of other agents.
  We take the point of view of a planner who optimizes the average profit of all agents, hence we fall into the family of the so-called ``Mean Field Control'' problems.
  The simultaneous presence of the carryover effect makes the problem infinite dimensional
  hence belonging to a family of problems which are very difficult in general and whose study started only very recently, see Cosso et Al, 2023 (\cite{Cosso-et-al2022}).
Here we consider, as a first step, a simple version of the problem providing the solutions in a simple case through a suitable auxiliary problem.
\end{abstract}

\textbf{Key words}:
Mean field control problems;
Optimal advertising models;
Delay in the control;
Infinite dimensional reformulation.


\tableofcontents

\section{Introduction}
Since the seminal papers of \cite{VW} and \cite{NA} on dynamics model in marketing, a considerable amount of work has been devoted to problems of optimal advertising,
 both in
 monopolistic and
 competitive settings, and
 both in
 deterministic and stochastic environments (see \cite{sethi} for a review of the existing work until the
 {1990's}).

Various extensions of the basic setting of \cite{VW} and \cite{NA} have been studied. For the stochastic case, we recall, among the various papers on the subject, \cite{sethi-VWcomp}, \cite{grosset}, \cite{gw}, \cite{MoPh}.

Our purpose here is to start exploring a family of models
{that put} together two important features that may arise in such problems and that have not yet been {satisfactorily treated} in the actual theory on optimal control.
\\
On one side we account, as in \cite{GM} and \cite{GMSJOTA}
for the presence of delay effects, in particular the fact that the advertisement spending affects the
goodwill with some delay, the so-called carryover effect (see e.g. \cite{hartl}, \cite{sethi}, \cite{GMSJOTA} and the references therein).

{On the other side, and more {crucially}, we take into account the fact that the agents maximizing their profit/utility from advertising are embedded in an environment where other agents act and where the action of such other agents influences their own outcome (see e.g. \cite{MoPh} for a specific case of such a situation).
To model such interaction among maximizing agents, one typically resorts to game theory. However, cases like this, where the number of agents can be quite large (in particular if we {hink of} web advertising), are very difficult to treat in an $N$-agents game {setting}.
A way to make such a problem tractable but still meaningful {is} to resort to what is called the mean-field theory. The idea is the following: assume that the agents are homogeneous (i.e.\ displaying the same state equations and the same objective functionals) and send their numbers to infinity.
The resulting limit problem is in general more treatable, and, under certain
conditions, its equilibria are a good approximation
of the
$N$-agents game (see e.g. the books \cite{CarmonaDelarueBook} for an extensive survey on the topic).
\\
For the above reason, we think it is interesting, both from the mathematical and economic side, to consider the optimal advertising investment problem with delay of \cite{GM,GMSJOTA} in the case when, in the state equation and in the objective, one adds a mean field term depending on the law of the state variable (the goodwill), which takes into account the presence of other agents.
\\
There are two main ways of looking at the problem when such mean field terms are present. One (which {falls} into the class of Mean Field Games (MFG), see e.g. \cite[Ch.\ 1]{CarmonaDelarueBook}, and which is not our goal here) is to look at the Nash equilibria where each agent takes the distribution of the state variables of the others as given. The other one, which we follow here, is to assume a cooperative game point of view: there is a planner that optimizes the average profit of each agent: this means that we fall into the family of the so-called ``Mean Field Control'' (MFC) problems (or ``control of McKean-Vlasov dynamics'').
We believe that both viewpoints are interesting from the economic side and challenging from the mathematical side.
In particular, the one we adopt here (the Mean Field Control) can be seen as a benchmark (a first best) to compare, subsequently, with the non-cooperative Mean Field Game case, as is typically done in game theory (see e.g. \cite{BFFGGEB}). It can also be seen as the case of {a} big selling company (who {acts} as the central planner), which has many shops in the territory whose local advertising policies interact.}

The simultaneous presence of the carryover effect and of the ``Mean Field Control'' terms makes the problem {belong} to the family of infinite dimensional control of McKean-Vlasov dynamics: a family of problems that are very difficult in general and whose study started only very recently (see \cite{Cosso-et-al2022}).

{Here we consider, as a first step, a simple version of the problem that displays
  {a}
  linear state equation, mean field terms depending only on the first moments, and an objective functional whose integrand (the running objective) is separated in the state and the control.
We develop the infinite dimensional setting in this case.
Moreover, we show that,
in the special subcase when the running objective is linear in the state and quadratic in the control, we can solve the problem.
This is done through the study of a suitable auxiliary problem whose HJB equation can be explicitly solved (see Section 4 below) and whose optimal feedback control can be found through an infinite dimensional Verification Theorem (see Section 5 below).}

The paper is organized as follows.
\begin{itemize}
  \item
In Section \ref{sec:formul}, we
formulate the optimal advertising problem as an optimal control
problem for 
stochastic delay differential equations with mean field terms and delay in the control.
Moreover, using that the mean field terms depends
only on the first moments we introduce an auxiliary problem without mean field terms but with a ``mean'' constraint on the control (see \eqref{eq:V_boh}).
\item
In Section \ref{sec:reformul},
the above ``not mean field''  auxiliary non-Markovian optimization problem is  ``lifted'' to an infinite dimensional Markovian control
problem, still with a ``mean'' constraint on the control (see \eqref{equivVxVY}).
\item
In Section \ref{sec:HJB}, we show how to solve the original problem in the special case when the optimal controls of the original and auxiliary {problems} are deterministic. We explain the strategy in Subsection \ref{SS:strategy}, proving Proposition \ref{pr:newFG}.
Then we consider a suitable Linear Quadratic (LQ) case. In Subsection \ref{sec:LQnew}, we solve the
  appropriated HJB equation,
  while, in Subsection \ref{sec:ver}, we find, through a verification theorem, the solution of the auxiliary LQ problem. Finally, in Subsection \ref{SS:back}, we show that we can use Proposition \ref{pr:newFG})
  to {also get} the solution of the original LQ problem.
\end{itemize}


%

\section{Formulation of the problem}\label{sec:formul}
{We call $X(t)$ the
stock of advertising goodwill (at time $t \in [0,T]$) of a given product.
We assume that the dynamics of $X(\cdot)$ is given by the following} controlled stochastic delay differential equation (SDDE), where $u$ models the intensity of advertising spending:
  \begin{equation}
\label{eq:SDDE}
\begin{dcases}
  dX(t) = \left[a_0 X(t) +a_1\E[ X(t)]
        + b_0 u(t) +
            \int_{-d}^0b_1(\xi)u(t+\xi) d\xi\right]dt + \sigma dW(t)& \forall t\in[0,T] \\[10pt]
X(0)=x \\[10pt]
u(\xi)=\delta(\xi)\qquad\forall\xi\in[-d,0]&
\end{dcases}
\end{equation}
\noindent
where the Brownian motion $W$ is defined on a filtered probability space $(\Omega,\mathcal{F},\mathbb{F}=(\mathcal{F}_t)_{t\geq 0},%
\mathbb{P})$,
  with $(\Omega,\mathcal{F},\mathbb{P})$ being complete,
$\mathbb{F}$ being the augmentation of the filtration generated by $W$,
and where, for a given closed interval $U\subset \mathbb{R}$, the control strategy
$u$ belongs to $\mathcal{U}\coloneqq L^2_\mathcal{P}(\Omega\times[0,T];U)$, the space of
$U$-valued square
integrable progressively measurable processes. 
The last line in~\eqref{eq:SDDE}
must read as an extension of $u$ to $[-d,T]$ by means of $\delta$.


{Here the control space and the state space are both equal to the set $\R$ of real numbers\footnote{This means that, due to the difficulty of the problem, we do not consider ex ante state or control constraints. They could be checked ex post or could be the subject of a subsequent research work.}}
{Regarding the coefficients and the initial data, we assume} the following conditions are verified:
\begin{assumption}\label{ip-concrete} $\quad$
\begin{enumerate}
\item $a_0,a_1\in \R$;
\item $b_0 \geq 0$;
\item $b_1(\cdot) \in L^2([-d,0];\R^+)$;
\item {$\delta(\cdot)\in L^2([-\delta,0];U)$}
\end{enumerate}
\end{assumption}
Here $a_0$ and $a_1$ are constant factors reflecting the goodwill changes in absence of advertising, $b_0$ is a constant advertising effectiveness factor,
and $b_1(\cdot)$
is the density function of the time lag between the
advertising expenditure $u$
and the corresponding effect on the
goodwill level. Moreover, $x$ is the level of goodwill at the
beginning of the advertising campaign, $\delta(\cdot)$ is the history of the advertising expenditure before time zero (one can
assume $\delta(\cdot)=0$, for instance).

Notice that under Assumption~\ref{ip-concrete} there exists a unique strong solution to
  the following SDDE starting at time $t\in[0,T)$:
\begin{equation}
\label{eq:SDDE-t}
\hspace{-0.4cm}\begin{dcases}
  dX(s) = \left[a_0 X(s) +a_1\E[ X(s)]
        + b_0 u(s) +
            \int_{-d}^0b_1(\xi)u(s+\xi) d\xi\right]ds + \sigma dW(s)& \forall s\in[t,T] \\[10pt]
X(t)=x\\[10pt]
u({t+\xi})=\delta(\xi)\quad\forall\xi\in[-d,0]&
\end{dcases}
\end{equation}
We denote such a solution by $X^{t,x,u}$.
It belongs $ L^2_\mathcal{P}(\Omega\times [0,T],\mathbb{R})$. In what follows, without loss of generality, we always assume to deal with a continuous version  $X^{t,x,u}$.

The objective functional to be maximized is defined as
\begin{equation}
  \label{eq:obj-orig}
  \begin{split}
    J(t,x;u(\cdot)) =& \E \left[ \int_t^T e^{-r(s-t)} \left( f
        \left(s, X^{t,x,u}(s),\mathbb{E} \left[ X^{t,x,u}(s) \right],u(s),\mathbb{E} \left[ u(s) \right]
          \right)
      \right)ds\right.\\
    &+e^{-r(T-t)} \mathbb{E}\left[ g \left( X^{t,x,u}(T),
        \mathbb{E}\left[ X^{t,x,u}(T) \right] \right) \right]
  \end{split}
\end{equation}
where for the functions
$ f\colon[0,T]{\times}\mathbb{R}\times \mathbb{R} \rightarrow
\mathbb{R}$ and
$ g\colon\mathbb{R}\times \mathbb{R} \rightarrow \mathbb{R}$ we assume
the following Assumption~\ref{ass-fg} is verified.
\begin{assumption}\label{ass-fg}$\quad$
  \begin{enumerate}
  \item The functions $f,g$
    are measurable.
  \item There exist $N>0,{\ell}>0, \theta>1$ such that
    \begin{equation*}
      f(t,x,m ,u,z )
      +
      g(x,m )
      \leq
      N(1+|x|+|m |+|u|+|z |)-{\ell}(|u|+|z |)^\theta,
    \end{equation*}
    for all $t\in[0,T],y\in \mathbb{R},m \in \mathbb{R},z \in \mathbb{R}$.
  \item $f,g$ are locally uniformly continuous in $x,m $, uniformly with respect to $(t,u ,z )$, meaning that for every $R>0$ there exists a modulus of continuity
    $\mathtt{w}_R\colon \mathbb{R}^+\rightarrow \mathbb{R}^+$ such that
    \begin{equation*}
      \sup_{\substack{t\in[0,T]\\u\in \mathbb{R},z \in \mathbb{R}}}|f(t,x,m ,u,z )-f(t,x',m ',u,z )|+|g(x,m )-g(x',m ')|
     \leq \mathtt{w}_R (|x-x'|+|m -m '|  )
   \end{equation*}
   for all real numbers $x,m ,x',m '$ such that $|x|\vee|m |\vee|x'|\vee|m '|\leq R$.
  \end{enumerate}
\end{assumption}
Under Assumption~\ref{ip-concrete} and
Assumption~\ref{ass-fg}, the reward functional $J$ in
\eqref{eq:obj-orig} is
well-defined for any $(t,x;u(\cdot))\in [0,T]\times \mathbb{R}^+\times
\mathcal{U}$.
\noindent We also define the value function $ \overline V$ for this problem as follows:
\begin{equation}\label{bar-value}
  \overline V(t,x) = \sup_{u  \in\mathcal{U}} J(t,x;u  ),
\end{equation}
{for $(t,x)\in[0,T]\times \mathbb{R}$.}
We shall say that $u^*  \in\mathcal{U}$ is an optimal control strategy if it is such that
$$ \overline V(t,x)=J(t,x;u^*  ). $$
Our main aim here is to finding such optimal control strategies

\bigskip

We now take into account the controlled ordinary delay differential equation (ODDE)
  \begin{equation}
\label{eq:ODDE-t}
\begin{dcases}
  dM(s) =  \left( (a_0+a_1) M(s)
        + b_0 z(s) +
              \int_{-d}^0b_1(\xi)z(s+\xi) d\xi \right)  ds &\qquad \forall s\in[t,T] \\[10pt]
M(t)=m\\[10pt]
z({t+\xi})=\delta(\xi)\qquad\forall\xi\in[-d,0]&
\end{dcases}
\end{equation}
where $m\in \mathbb{R}$ and $z\in L^2([0,T],\mathbb{R})$ is extended to $[-d,0]$ by $\delta$ as expressed by the last line in~\eqref{eq:ODDE-t}.
We denote by $M^{t,m,z}$ the unique strong solution to \eqref{eq:ODDE-t}.
It is straightforward to notice the relationship
\begin{equation}
  \label{eq:XM-rel}
  M^{t,m,z}=\mathbb{E} \left[ X^{t,m,u} \right]
  \
  \mathrm{whenever}
  \ 
  z(s)=\mathbb{E}[u(s)]\ \mathrm{for}\ s\in[t,T].
\end{equation}
Property \eqref{eq:XM-rel} suggests that we can couple
the two systems
\eqref{eq:SDDE-t}
and \eqref{eq:ODDE-t}
as follows.
We set
\begin{equation}\label{A0}
  A_0\coloneqq      \begin{bmatrix}
      a_0&a_1\\
      0&a_0+a_1
    \end{bmatrix}
  \end{equation}
  and introduce, for $\tilde{x}\in \mathbb{R}^2$ and with
  \begin{equation}\label{tilde-usigma}
    \tilde{u}=(u,z)\in \tilde{\mathcal{U}}\coloneqq
    L^2_\mathcal{P} \left( \Omega\times [0,T];\mathbb{R} \right) \times L^2 \left([0,T];\mathbb{R}  \right),\qquad
      \tilde{\sigma}= (\sigma,0),
    \end{equation}
    the process $\tilde{X}^{t,\tilde{x},\tilde{u}}$ as the unique strong solution of the
controlled SDDE
\begin{equation}
  \label{eq:SDDE-m}
 \hspace{-0.6cm} \begin{dcases}
    d\tilde{X}(s)=
     \left(
   A_0    \tilde{X}(s)+b_0 \tilde{u}(s)
       +
       \int_{-d}^0
       b_1(\xi)
       \tilde{u}(s+\xi)
       d\xi
     \right) ds
     +
     \tilde{\sigma}
     dW(s)& \forall s\in(t,T]\\[10pt]
     \tilde{X}(t)=\tilde{x}
     &\\[10pt]
     \tilde{u}({t+\xi})= \left( \delta(\xi),\delta(\xi) \right)
     \quad \forall \xi\in[-d,0]&
  \end{dcases}
\end{equation}
then
by \eqref{eq:SDDE-t},
\eqref{eq:ODDE-t}
,
\eqref{eq:XM-rel},
and \eqref{eq:SDDE-m},
we immediately have
\begin{equation}
  \label{eq:Xtilde-XM}
 \left(       X^{t,x,u},M^{t,x,z} \right)
=
    \tilde{X}^{t,(x,x),\tilde{u})}
  \  \underline{\mathrm{if}
  \ z(s)=\mathbb{E}[u(s)]\ \mathrm{for}\ s\in[t,T]}.
\end{equation}
Property \eqref{eq:Xtilde-XM}
states that the process $X^{t,x,u}$ can be seen as the first projection of a bidimensional process
driven by a SDDE whose coefficients do not involve any dependence on the law.

Thanks to \eqref{eq:Xtilde-XM}, we can rephrase the original control problem as follows.
We define, for $t\in[0,T],\tilde{x}\in \mathbb{R}^2$, and for
\begin{equation*}
  \tilde{u}=(u,z)\in \tilde{\mathcal{U}}\coloneqq
  L^2_\mathcal{P} \left( \Omega\times [0,T];\mathbb{R} \right) \times L^2 \left([0,T];\mathbb{R}  \right),
\end{equation*}
the functional
  \begin{equation}
  \label{eq:obj-orig-tilde}
  \begin{split}
    \tilde{J}(t,\tilde{x};\tilde{u})\coloneqq &
    \E \left[ \int_t^T
      e^{-r(s-t)}
        f \left( s, \tilde{X}^{t,\tilde{x},\tilde{u}}(s)
        ,\tilde{u}(s) \right)
      ds+
      g \left( \tilde{X}^{t,\tilde{x},\tilde{u}}(T) \right)  \right],
  \end{split}
\end{equation}
where, with a slight abuse of notation, we identify
\begin{equation}\label{fg_notation-couples}
  f(t,(x,m),(u,z))=f(t,x,m,u,z)\qquad g((x,m))=g(x,m).
\end{equation}

Then, by
\eqref{eq:obj-orig},
\eqref{bar-value},
\eqref{eq:Xtilde-XM}, and
\eqref{eq:obj-orig-tilde},
it follows that
\begin{equation}
  \label{eq:V_boh}
  \overline{V}(t,x)=\sup
  \left\{
    \tilde{J}(t,(x,x);\tilde{u})\colon
    {
      \tilde{u}\in \tilde{\mathcal{U}},
    \ \mathrm{and}\ }
    z(s)=\mathbb{E}[u(s)]\ s\in[t,T]
  \right\}.
\end{equation}


\section{Carryover effect of advertising: reformulation of the {problem in infinite dimension}}
\label{sec:reformul}


%
To recast the
{SDDE (\ref{eq:SDDE-m})} as an abstract stochastic differential equation
on a suitable Hilbert space we use the approach introduced first by \cite{VK} in the deterministic case and then
extended in \cite{GM} to the stochastic case (see also \cite{FGFM-III} where the case of unbounded control operator is considered).
We reformulate equation (\ref{eq:SDDE-m}) as an abstract stochastic differential equation in the following Hilbert space $H$
\begin{equation*}
  H\coloneqq \R^2\times L^2([-d,0],\R^2).
\end{equation*}
If $y\in H$, we denote by $y_0$ the projection of $y$ onto $\mathbb{R}^2$ and by $y_1$
the projection of $y$ onto $L^2([-d,0],\mathbb{R}^2)$. Hence $y=(y_0,y_1)$.
The inner product in $H$ is induced by its factors, meaning
$$ \langle y,y'\rangle \coloneqq \<y_0,y'_0\>_{\R^2} + \int_{-d}^0 \<y_1(\xi),y_1'(\xi)\>_{\R^2} d\xi\qquad \forall y,y'\in H.
$$
In particular, the induced norm is
$$
|y| = \left( |y_0|_{\R^2}^2 + \int_{-d}^0 |y_1(\xi)|_{\R^2}^2 d\xi\right)^{1/2}\qquad \forall y\in H.
$$


\noindent
Recalling~\eqref{A0},
we define
$A\colon \mathcal{D}(A)\subset H\rightarrow H$ by
\begin{equation*}
{  Ay\coloneqq  \left(
A_0y_0,-\dot{y}_1
  \right)}
\end{equation*}
where the domain $\mathcal{D}(A)$ is
\begin{equation*}
  \mathcal{D}(A)=
  \left\{
    y\in H\colon y_1\in W^{1,2}([-d,0],\mathbb{R}^2),
    \ y_1(-d)=0
  \right\}.
\end{equation*}
The adjoint $A^*\colon \mathcal{D}(A^*)\subset H\rightarrow H$ of $A$  is given by
\begin{equation*}
{ A^*y\coloneqq  \left(
A_0^*y_0,\dot{y}_1
  \right)}
\end{equation*}
with
\begin{equation*}
  \mathcal{D}(A^*)=
  \left\{
    y\in H\colon y_1\in W^{1,2}([-d,0],\mathbb{R}^2),
    \ y_1(0)=y_0
  \right\}.
\end{equation*}
The operator $A$ generates a $C_0$-semigroup $\{e^{tA}\}_{t\in \mathbb{R}^+}$ on $H$, where
\begin{equation*}
  e^{tA}y=
  \left(
    e^{tA_0}y_0+\int_{-d}^0 \mathbf{1}_{[-t,0]}e^{(t+s)A_0}y_1(s)ds,
    y_1(\cdot-t)\mathbf{1}_{[-d+t,0]}(\cdot)
  \right)
  \qquad \forall y\in H,
\end{equation*}
  whereas the $C_0$-semigroup
  $\{e^{tA^*}\}_{t\in \mathbb{R}^+}$ generated by $A^*$ is given by
  \begin{equation*}
    e^{tA^*}y=
    \left(
      e^{tA^*_0}y_0,
      e^{(\cdot+t)A^*_0}y_0\mathbf{1}_{[-t,0]}(\cdot)
      +
      y_1(\cdot+t)\mathbf{1}_{[-d,-t]}(\cdot)
    \right)
    \qquad \forall y\in H,
  \end{equation*}
  where $A^*_0$ is the adjoint of $A_0$.

  We then introduce the noise operator
  $G\colon \mathbb{R}\rightarrow H$
  defined by
  \begin{equation*}
    Gx\coloneqq
       \left( (\sigma x,        0),0    \right) \qquad \forall x\in \mathbb{R},
  \end{equation*}
  and the
  control operator
  $B\colon \mathbb{R}^2\rightarrow H$ defined by
  \begin{equation*}
    By_0=(b_0y_0,b_1(\cdot)y_0)\qquad \forall y_0\in \mathbb{R}^2.
  \end{equation*}
The adjoint $B^*\colon H\rightarrow \mathbb{R}^2$ of $B$ is given by
\begin{equation*}
  B^*y= b_0y_0+\int_{-d}^0b_1(\xi)y_1(\xi)d\xi\qquad \forall y\in H.
\end{equation*}
We now introduce the abstract stochastic differential equation on $H$
\begin{equation}\label{eq:abstract}
    \begin{dcases}
      dY(s)= \left( AY(s)+B \tilde{u}(s) \right)   ds+GdW(s)& s\in(t,T]\\[10pt]
      Y(t)=y&\\[10pt]
      \tilde{u}(t+\xi)= \left( \delta(\xi),\delta(\xi) \right)
     \quad \forall \xi\in[-d,0]&
  \end{dcases}
\end{equation}
with $t\in[0,T), y\in H,
\tilde{u}\in\mathcal{U}\times \mathcal{U}$.
Denote by $Y^{t,y,\tilde{u}}$ the mild solution to~\eqref{eq:abstract}, i.e., the pathwise continuous process in $L^2_\mathcal{P}(\Omega\times [0,T];H)$ given by the variation of constants formula:
\begin{equation}\label{eq-mild}
Y^{t,y,\tilde{u}}(s) = e^{(s-t)A}y + \int_t^s e^{(s-r)A}B\tilde{u}(r) dr+ \int_t^s e^{(s-t)A}G dW(r),\qquad  \forall s\in[t,T].
\end{equation}
Similarly as done in~\cite{GM},
if the space of admissible controls is restricted to
$ \tilde{\mathcal{U}}$,
one can show that
\eqref{eq:abstract} is equivalent to
\eqref{eq:SDDE-m}, in the sense that
\begin{equation}
  \label{eq:equiv-tildeXY}
  Y^{t,y,\tilde{u}}_0(s)=
  \tilde{X}^{t,y_0,\tilde{u}}
\end{equation}
for every $t\in[0,T),\tilde{u}\in \tilde{\mathcal{U}}$, and for every $y=(y_0,y_1)\in H$ with
\begin{equation}\label{y1-starting}
  y_1(\xi)=
  \left(
    \int_{-d}^\xi b_1(\zeta)\delta(\zeta-\xi)d\zeta,
    \int_{-d}^\xi b_1(\zeta)\delta(\zeta-\xi)d\zeta
  \right) \qquad \forall \xi\in[-d,0].
\end{equation}
A further equivalence is given by considering together
\eqref{eq:Xtilde-XM} and
\eqref{y1-starting}, that provide
\begin{equation}
  \label{eq:equiv-YXM}
   Y^{t,y,\tilde{u}}_0(s)=
 \left(       X^{t,x,u},M^{t,x,z} \right)
  \  \mathrm{if}
  \ y_0=(x,x),\ y_1\ \textrm{is as in }\eqref{y1-starting},\ z(s)=\mathbb{E}[u(s)]\ \mathrm{for}\ s\in[t,T].
\end{equation}
Thanks to equivalence
\eqref{eq:equiv-YXM}, we can rephrase the original control problem as follows.
For $t\in[0,T],y\in H, \tilde{u}\in \mathcal{U}\times\mathcal{U}$, define
the functional (recall~\eqref{fg_notation-couples})
  \begin{equation}
  \label{eq:obj-orig-Y}
  \begin{split}
    \mathcal{J}(t,y;\tilde{u})\coloneqq &
    \E \left[ \int_t^T
      e^{-r(s-t)}
        f \left( s, Y^{t,y,\tilde{u}}_0(s)
        ,\tilde{u}(s) \right)
      ds+
      g \left( Y^{t,y,\tilde{u}}_0(T) \right)  \right]
  \end{split}
\end{equation}
Then, by
\eqref{eq:obj-orig-tilde},
\eqref{eq:V_boh},
\eqref{eq:equiv-tildeXY},
and
\eqref{y1-starting},
it follows that
\begin{equation}\label{equivVxVY}
  \overline{V}(t,x)=
  \sup
  \big\{
    \mathcal{J}(t,y;\tilde{u})\colon y_0=(x,x),
    \ y_1\ \textrm{is as in }\eqref{y1-starting},\,
    \tilde{u}\in \tilde{\mathcal{U}},\,
    \mathrm{and}\ z(s)=\mathbb{E}[u(s)]\ s\in[t,T]
  \big\}.
\end{equation}

\section{Solution of the original problem in a special Linear Quadratic (LQ) case}
\label{sec:HJB}


\subsection{The strategy of solution through a suitable HJB equation}
\label{SS:strategy}

Following \eqref{equivVxVY} above we introduce the function
\begin{equation*}
  \mathcal{V}\colon [0,T]\times H\rightarrow \mathbb{R}
\end{equation*}
defined by
\begin{equation*}
\mathcal{V}(t,y)\coloneqq    \sup
  \big\{
    \mathcal{J}(t,y,\tilde{u})\colon\, \tilde{u}\in \tilde{\mathcal{U}},\;  z(s)=\mathbb{E}[u(s)]\ \forall s\in[t,T]
  \big\}.
\end{equation*}
Notice that, by
\eqref{equivVxVY}, we have
\begin{equation}
  \label{eq:equivVcalV}
  \underline{
    \overline{V}(t,x)
    =
    \mathcal{V}(t,y)
    \ \mathrm{if}
    \ y_0=(x,x),
    \ \mathrm{and\ if}
    \ y_1\ \textrm{is as in }\eqref{y1-starting}.
  }
\end{equation}

    The problem with the above constraint $z(s)=\mathbb{E}[u(s)]$, for $s\in[t,T]$, is that it does not allow to apply directly the Dynamic Programming Approach
to get the HJB equation.
For this reason,
instead of optimizing on the set $\mathcal{U}$ with the constraints
$z(s)=\mathbb{E}[u(s)]\ s\in[t,T]$, we
take into consideration a different problem, for which the optimization is performed on the set
$\mathcal{U}\times\mathcal{U}$ with the
constraint
$z(s)=u(s)\ s\in[t,T]$, hence considering the
following value function
\begin{equation}\label{optprV}
  V(t,y)
\coloneqq    \sup
  \big\{
  \mathcal{J}(t,y,\tilde{u})\colon \, \tilde{u}=(u,z)\in \mathcal{U}\times \mathcal{U},\
     \mathrm{and}\ u=z
  \big\}.
\end{equation}
In general we do not know if and how this function is related to
$\mathcal{V}$ (and consequently to our goal $\overline{V}$). However it is clear from the constraints involved that,  if
for both problems  $V$ and $\mathcal{V}$
the supremum is reached on the set of deterministic controls, meaning
\begin{subequations}
  \begin{align}
      \mathcal{V}(t,y)=&\mathrm{(to\ prove)}= \sup \big\{
      \mathcal{J}(t,y,\tilde{u})\colon \, \tilde{u}=(u,z)\in
      \mathcal{U}\times \mathcal{U},\ \mathrm{and}\ u=z\
      \mathrm{deterministic}
      \big\}\label{supDiag-a}\\[10pt]
      V(t,y)=&\mathrm{(to\ prove)}= \sup \big\{
      \mathcal{J}(t,y,\tilde{u})\colon \, \tilde{u}=(u,z)\in
      \mathcal{U}\times \mathcal{U},\ \mathrm{and}\ u=z\
      \mathrm{deterministic} \big\},\label{supDiag-b}
  \end{align}
\end{subequations}
then  finding the deterministic optimal controls for $\mathcal{V}$ is equivalent to doing that for ${V}$.
For future reference, we restate this observation in the following proposition.

\begin{proposition}\label{pr:newFG}
  Let $t\in [0,T]$ and $y\in H$.  If \eqref{supDiag-a} and
  \eqref{supDiag-b} hold true, then a deterministic control
  $\tilde{u}^*=(u^*,u^*)\in \mathcal{U}\times \mathcal{U}$ is optimal for
  $\mathcal{V}$ if and only if it is optimal for $V$.
\end{proposition}

\noindent The HJB equation associated to the optimal control problem related to ${V}$ is the following.
\begin{equation}
  \label{eq:HJVcalVH0}
  \begin{dcases}
     v _t(t,y)+\frac{1}{2}\tr Q\nabla^2 v (t,y)+\langle Ay,\nabla v (t,y) \rangle&\\
    \hspace{100pt}+H_0(t,y,\nabla  v (t,y))
    -r  v (t,y)
    =0&
     \forall (t,y)\in (0,T)\times H
    \\[10pt]
     v (T,y)=g(y_0) \qquad  \forall y\in H&
  \end{dcases}
\end{equation}
where $Q=G^*G$,
and the Hamiltonian function defined as
\begin{equation*}
  H_0(t,y,p)
  \coloneqq
  \sup_{\tilde{u}\in \mathbf{D}}
  H_{CV}(t,y,\tilde{u},p)=
    \sup_{\tilde{u}\in \mathbf{D}}
  \big\{
    f(t,y_0,\tilde{u})+\langle B\tilde{u},p\rangle
  \big\},
\end{equation*}
with $H_{CV}$ denoting the current value Hamiltonian function,
and $\mathbf{D}$ being the diagonal in $U\times U$, meaning
  $\mathbf{D}= \left\{ (u,u)\colon u\in U \right\}$.
Notice that $H_0(t,y,p)$ depends on $p$ only by means of $B^*p$. Indeed, if we define
\begin{equation}\label{HamFunc}
  H(t,y,q)
  \coloneqq
    \sup_{\tilde{u}\in \mathbf{D}}
  \big\{
    f(t,y_0,\tilde{u})+\langle \tilde{u},q\rangle
  \big\},
\end{equation}
we get $H_0(t,y,p)=H(t,y,B^*p)$.
Then \eqref{eq:HJVcalVH0} can be rewritten as
\begin{equation}
  \label{eq:HJVcalVH}
  \begin{dcases}
     v _t(t,y)+\frac{1}{2}\tr Q\nabla^2 v (t,y)+\langle Ay,\nabla v (t,y) \rangle&\\
    \hspace{100pt}+H(t,y,B^*\nabla  v (t,y))
    -r  v (t,y)
    =0&
     \forall (t,y)\in (0,T)\times H
    \\[10pt]
     v (T,y)=g(y_0) \qquad  \forall y\in H&
  \end{dcases}
\end{equation}
  Notice that, in the above equations
  \eqref{eq:HJVcalVH0} and \eqref{eq:HJVcalVH},
  the gradient inside the Hamiltonian $H$ is indeed a couple of directional derivatives since it acts only through the operator $B^*$ whose image lies in $\mathbb{R}^2$.

In the next subsections we specify $f,g$ and we show that with such a choice
\eqref{supDiag-a}
and
\eqref{supDiag-b} are verified.



\subsection{Explicit solution of the HJB equation in the auxiliary LQ case}
\label{sec:LQnew}


  In this section we specify
  the general model with
  \begin{equation}\label{lq-case}
    \begin{split}
      f(t,x,m ,u,z ) &= \alpha_0 x-\alpha_1 m -\beta_0u-\gamma_0 u^2-\beta_1 z -\gamma_1z^2\\
      g(x,m ) &=  \lambda_0 x-\lambda_1 m
    \end{split}
  \end{equation}
  for $(x,m,u,z)\in \mathbb{R}^4$,
  where
  \begin{enumerate}
  \item $\alpha_0,\alpha_1,\beta_0,\beta_1,\lambda_0, \lambda_1\in \mathbb{R}$;
  \item $\gamma_0>0, \gamma_1>0$.
  \end{enumerate}
  \noindent We also set $U=\mathbb{R}$. Notice that
  Assumption~\ref{ass-fg} is satisfied.
  Moreover, denoting $\tilde{\alpha}=(\alpha_0,-\alpha_1)$,
  $\tilde{\beta}=(\beta_0,\beta_1)$,
  and recalling~\eqref{fg_notation-couples}, we have, for $q\in \mathbb{R}^2$,
  \begin{equation}
    \label{ustar}
    \begin{split}
      u^*(q)\coloneqq &
      \argmax_{u\in U}
      \big \{
      \langle \tilde{\alpha}, y_0\rangle
      +\langle q-\tilde{\beta},(1,1)\rangle u
      -(\gamma_0+\gamma_1)u^2
      \big\}\\
      =&\frac{\langle q-\tilde{\beta},(1,1)\rangle}{2(\gamma_0+\gamma_1)},
    \end{split}
  \end{equation}
  {which} entails, by considering the definition of $H$ given in \eqref{HamFunc},
  \begin{equation*}
    H(t,y,q)=\frac{ \left( \langle q-\tilde{\beta},(1,1)\rangle \right) ^2}{4(\gamma_0+\gamma_1)}+\langle\tilde{\alpha},y_0\rangle
  \end{equation*}
  and then the HJB equation
\eqref{eq:HJVcalVH0} reads as
\begin{equation}
  \label{eq:HJVcalVH_particolare}
  \begin{dcases}
     v _t(t,y)+\frac{1}{2}\tr Q\nabla^2 v (t,y)+\langle Ay,\nabla v (t,y) \rangle&\\
    \hspace{30pt}+
\frac{ \left( \langle B^*\nabla  v (t,y)-\tilde{\beta},(1,1)\rangle \right) ^2}{4(\gamma_0+\gamma_1)}+\langle\tilde{\alpha},y_0\rangle
    -r  v (t,y)
    =0&
     \forall (t,y)\in (0,T)\times H
    \\[10pt]
     v (T,y)=\langle \tilde{\lambda},y_0\rangle \qquad  \forall y\in H&
  \end{dcases}
\end{equation}
where $\tilde{\lambda}=(\lambda_0,-\lambda_1)$.
\noindent We look for solutions of~\eqref{eq:HJVcalVH_particolare}
of the following form
\begin{equation}\label{sol-forma-special}
   v (t,y)=\langle a(t),y\rangle+b(t)
\end{equation}
with $a\colon [0,T]\rightarrow H$ and $b\colon [0,T]\rightarrow \mathbb{R}$ to be determined.
The final condition in~\eqref{eq:HJVcalVH_particolare} holds true for~\eqref{sol-forma-special} only if
\begin{equation}
  \label{eq:aT_bT}
  a(T)=(\tilde{\lambda},0),\quad b(T)=0.
\end{equation}
Moreover, if $ v $ is of the form~\eqref{sol-forma-special},
\eqref{eq:HJVcalVH_particolare} reads as
\begin{equation}\label{eq:HJB-special}
  \langle \dot a(t),y\rangle+\dot b(t)+\langle y, A^*a(t)\rangle
  +
  \frac{ \left( \langle B^*a(t)-\tilde{\beta},(1,1)\rangle \right) ^2}{4(\gamma_0+\gamma_1)}+\langle\tilde{\alpha},y_0\rangle
  -r\langle a(t),y\rangle-rb(t)=0
\end{equation}
The previous equation \eqref{eq:HJB-special} is to be intended in a mild way that we are going to specify in the following, since we cannot guarantee that, for all $t$, $a(t)\in\cald(A^*)$. Indeed, by 
\eqref{eq:aT_bT},
$a(T)\notin \cald(A^*)$.

Equation~\eqref{eq:HJB-special} can be {split} into two equations by {isolating} the terms containing $y$ and all the other terms, namely
\begin{equation}\label{eq:HJB-special_a}
  \langle \dot a(t),y\rangle
  +\langle y, A^*a(t)\rangle
  +\langle\tilde{\alpha},y_0\rangle
  -r\langle a(t),y\rangle=0
\end{equation}
and
\begin{equation}\label{eq:HJB-special_b}
  \dot b(t)
  +
  \frac{ \left( \langle B^*a(t)-\tilde{\beta},(1,1)\rangle \right) ^2}{4(\gamma_0+\gamma_1)}
  -rb(t)=0.
\end{equation}
Taking into account that~\eqref{eq:HJB-special_a} must hold for all $y\in H$, and combining~\eqref{eq:HJB-special_a} and \eqref{eq:HJB-special_b} with the final conditions \eqref{eq:aT_bT}, we obtain two separated equations, one for $a$ and one for $b$, namely
\begin{equation}
  \label{eq:HJB-a}
  \begin{dcases}
   \dot a(t)
  + A^*a(t)
  +(\tilde{\alpha},0)
  -r a(t)=0& t\in[0,T)\\[5pt]
    a(T)=(\tilde{\lambda},0)&
  \end{dcases}
\end{equation}
and
\begin{equation}
  \label{eq:HJB-b}
  \begin{dcases}
      \dot b(t)
  +
  \frac{ \left( \langle B^*a(t)-\tilde{\beta},(1,1)\rangle \right) ^2}{4(\gamma_0+\gamma_1)}
  -rb(t)=0& t\in[0,T)\\[5pt]
    b(T)=0&
  \end{dcases}
\end{equation}
\noindent We solve~\eqref{eq:HJB-a}, which turns out to be an abstract evolution equation in $H$, in mild sense, getting
\begin{equation}\label{eq:am-sol}
  a(t)  =e^{(T-t)(A^*-r)}(\tilde{\lambda},0)
  +\int_t^T
  e^{(s-t)(A^*-r)}
  (\tilde{\alpha},0)
  ds.
\end{equation}
Consequently we can write the solution to \eqref{eq:HJB-b}
\begin{equation}\label{eq:bm-sol}
  b(t) =\int_t^T \frac{1}{2}e^{-r(s-t)}
  \frac{ \left( \langle B^*a(s)-\tilde{\beta},(1,1)\rangle \right) ^2}{4(\gamma_0+\gamma_1)} ds,
\end{equation}
where $a$ is given by \eqref{eq:am-sol}.

So far we have found a solution $v$ to the HJB equation~\eqref{eq:HJVcalVH_particolare} whose candidate optimal feedback is deterministic. In the next section we will prove that
it is indeed the optimal control and that $v=V$.
We will also prove that the optimal feedback control associated to the optimal control problem associated to $\mathcal{V}$ is deterministic.
This will allow us to apply Proposition \ref{pr:newFG}, so finding the optimal strategies for the initial problem in the linear quadratic case.

\subsection{Fundamental Identity and Verification Theorem in the auxiliary LQ case}
\label{sec:ver}

  The aim of this subsection is to provide a verification theorem and the existence of optimal feedback controls for the linear quadratic problem for $V$ introduced in the previous section.
  This, in particular, will imply that the solution in (\ref{sol-forma-special}), with $a$ and $b$ given respectively by \eqref{eq:am-sol} and \eqref{eq:bm-sol},  coincides with the value function of our optimal control problem
  $V$ defined in   \eqref{optprV}.

  The main tool needed to get the wanted results is an identity (often called ``{\em fundamental identity}'', see equation  (\ref{relfond})) satisfied by the solutions of the HJB equation. Since the solution (\ref{sol-forma-special}) is not smooth enough
  (it is not differentiable with respect to $t$
  due to the presence of $A^*$ in $a$, given by \eqref{eq:am-sol}), we  need to perform an approximation procedure thanks to which Ito's formula can be applied.
Finally we pass to the limit and obtain the needed ``{\em fundamental identity}''.





\begin{proposition}\label{prop rel fond}
  Let Assumption~\ref{ip-concrete} hold.
  Let $v$ be  as in \eqref{sol-forma-special}, with $a$ and $b$ given respectively by \eqref{eq:am-sol} and \eqref{eq:bm-sol}, solution of the HJB equation\eqref{eq:HJVcalVH_particolare}.
Then for every $t\in[ 0,T],\, y\in H$, and
     $\tilde{u}=(u,z)\in \mathcal{U}\times \mathcal{U}$, with  $u=z$,
 we
 have the fundamental identity
\begin{equation}\label{relfond}
  \begin{split}
    v(t,y)
    =&
    \mathcal{J}(t,y;\tilde{u})
    +\mathbb{E}
     \left[ \int_t^T
       e^{-r(s-t)}
       \left(
         \frac{ \left( \langle B^*\nabla  v (t,Y ^{t,y,\tilde{u}}(s))-\tilde{\beta},(1,1)\rangle \right) ^2}{4(\gamma_0+\gamma_1)}+
          \right.\right.\\
      &\left.\left.
          \phantom{
            \frac{\left(\tilde{\beta}B\right)^2}{(\gamma_0)}
          }
       +\langle\tilde{\alpha},(Y ^{t,y,\tilde{u}})_0(s)\rangle   -
          H_{CV}(s,B^*\nabla v(s,Y^{t,y,\tilde{u}},\tilde{u}(s))
           \right) ds
         \right].
%
  \end{split}
\end{equation}
\end{proposition}
%

\dim
  Let $t\in [0,T), y\in H, \tilde{u}=(u,z)\in \mathcal{U}\times \mathcal{U}$, $u=z$.
We should apply Ito's formula to the process
$ \left\{ e^{-rs}v(s,Y^{t,y,\tilde{u}}(s)) \right\} _{s\in[t,T]}$, but we cannot, because $Y^{t,y,\tilde{u}}$ is a mild solution
 (the integrals in  \eqref{eq-mild} are convolutions with a $C_0$-semigroup) and not a strong solution of
 \eqref{eq:abstract}, moreover $v$ is not differentiable in $t$, since $(\tilde{\lambda},0)\not\in D(A^*)$.
 Then we approximate $Y^{t,y,\tilde{u}}$ by means of the Yosida approximation
(see also \cite{FGFM-II}[Proposition 5.1]).
For $k_0\in \mathbb{N}$ large enough, the operator $k-A$, $k\geq k_0$, is full-range and invertible, with continuous inverse, and $k(k-A)^{-1}A$ can be extended to a continuous operator on $H$. Define, for $k\geq k_0$, the operator on $H$
 \begin{equation*}
   A_k\coloneqq k(k-A)^{-1}A.
 \end{equation*}
 It is well known that, as $k\rightarrow \infty$, $e^{tA_k}y'\rightarrow e^{tA}y'$ in {$H$},  uniformly
 for  $t\in [0,T]$ and for $y'$ on compact sets of $H$.
Since $A_k$ is continuous, there exists a unique strong solution $Y^{t,y,\tilde{u}}_k$ to the SDE on $H$
\begin{equation}\label{eq:abstract-Yosida}
    \begin{dcases}
      dY_k(s)= \left( A_kY_k(s)+B \tilde{u}(s) \right)   ds+GdW(s)& s\in(t,T]\\[10pt]
      Y_k(t)=y&\\[10pt]
      \tilde{u}(s+\xi)= \left( \delta(\xi),\delta(\xi) \right)
     \quad \forall \xi\in[-d,0]&
  \end{dcases}
\end{equation}
By taking into account \eqref{eq-mild} together with the
same formula with $A_k$ in place of $A$, and by recalling the convergence $e^{\cdot A_k}\rightarrow e^{\cdot A}$ mentioned above,
one can easily show that
\begin{equation}
  \label{eq:conv-Y_k-Y}
  Y_k^{t,y,\tilde{u}}\rightarrow Y^{t,y,\tilde{u}}
  \
  \mathrm{in}
  \ L_\mathcal{P}^{2}(\Omega\times [0,T];H)
  \ \mathrm{as}\ k\rightarrow \infty.
\end{equation}

\noindent We now take into consideration the HJB
  \begin{equation}
  \label{eq:HJVcalVH_particolare-k}
 \begin{dcases}
    v _t(t,y)+\frac{1}{2}\tr Q\nabla^2 v (t,y)+\langle A_ky,\nabla v (t,y) \rangle&\\    \hspace{30pt}+
\frac{ \left( \langle B^*\nabla  v (t,y)-\tilde{\beta},(1,1)\rangle \right) ^2}{4(\gamma_0+\gamma_1)}+\langle\tilde{\alpha},y_0\rangle
    -r  v (t,y)
    =0&
     \forall (t,y)\in (0,T)\times H    \\[10pt]
     v (T,y)=\langle \tilde{\lambda},y_0\rangle \qquad  \forall y\in H.&
  \end{dcases}
\end{equation}
As argued for
\eqref{eq:HJVcalVH_particolare}, a solution for
\eqref{eq:HJVcalVH_particolare-k}
  is given by
\begin{equation}\label{sol-forma-special-k}
   v^{(k)} (t,y)=\langle a_k(t),y\rangle+b_k(t)
 \end{equation}
 where
 \begin{equation}\label{eq:am-sol-k}
  a_k(t)  =e^{(T-t)(A_k^*-r)}(\tilde{\lambda},0)
  +\int_t^T
  e^{(s-t)(A_k^*-r)}
  (\tilde{\alpha},0)
  ds
\end{equation}
and
\begin{equation}\label{eq:bm-sol-k}
  b_k(t) =\int_t^T \frac{1}{2}e^{-r(s-t)}
  \frac{ \left( \langle B^*a_k(s)-\tilde{\beta},(1,1)\rangle \right) ^2}{4(\gamma_0+\gamma_1)} ds.
\end{equation}
Since $A_k^*\in L(H)$, both $a_k$ and $b_k$ belong
to $C^{1}([0,T];\mathbb{R})$.
 So we can appy Ito's formula to
$ \left\{ e^{-r(s-t)}v^{(k)}(s, Y_k^{t,y,\tilde{u}}(s)) \right\} _{s\in[t,T]}$ getting:
\begin{equation*}
  \begin{split}
    e^{-r(T-t)}&\E  \left[ v^{(k)}(T,Y_k^{t,y,\tilde{u}}) \right]  -\E \left[  v^{(k)}(t,y) \right] \\
    =&\E
    \left[
    \int_t^Te^{-r(s-t)}\left( v^{(k)}_t(s,Y_k^{t,y,\tilde{u}}(s))-rv^{(k)}(s,Y_k^{t,y,\tilde{u}}(s))+ {1\over 2}\tr \left[ Q \nabla^2v^{(k)}(t,Y_k^{t,y,\tilde{u}}(s)) \right]  \right. \right.\\
  &\qquad\phantom{\frac{1}{2}}\left.+
      \langle
      A_kY_k^{t,y,\tilde{u}}(s),
    \nabla v^{(k)}
    (s,Y_k^{t,y,\tilde{u}}(s))
      \rangle
      +\langle B\tilde {u}(s),\nabla v^{(k)}(s,Y_k^{t,y,\tilde{u}}(s))\rangle\right]ds.
\end{split}
\end{equation*}
Since $v^{(k)}$ is a solution to  equation \eqref{eq:HJVcalVH_particolare-k}, we get
\begin{equation}\label{quasirelfondv^nk}
  \begin{split}
    e^{-r(T-t)}& \E \left[ \langle\tilde{\lambda},
      \left(Y_k^{t,y,\tilde{u}}\right)_0(T)\rangle \right] -
      v^{(k)}(t,y)
    \\
    &=\E
    \int_t^T  \left[ e^{-r(s-t)}\left(
      -\frac{ \left( \langle B^*\nabla  v^{(k)} (t,Y_k^{t,y,\tilde{u}}(s))-\tilde{\beta},(1,1)\rangle \right) ^2}{4(\gamma_0+\gamma_1)}-
      \langle\tilde{\alpha},(Y_k^{t,y,\tilde{u}})_0(s)\rangle\right.\right.\\
    &\left.\left.\phantom{\frac{ -\left( \langle B^*\nabla v^{k}
            (t,Y_k^{t,y,\tilde{u}}(s))-\tilde{\beta},(1,1)\rangle
          \right) ^2}{4(\gamma_0+\gamma_1)}} +\<B\tilde {u}(s),\nabla
      v^{(k)}(s,Y_k^{t,y,\tilde{u}}(s))\>\right)ds \right].
  \end{split}
\end{equation}
We then let $k\rightarrow\infty$ in
\eqref{quasirelfondv^nk}.
Recalling the
convergence $e^{\cdot A_k}\rightarrow e^{\cdot A}$ mentioned above, we
first notice that
\begin{equation}\label{ak-a_bk-b}
  a_k\rightarrow a
  \ \mathrm{in}\ H
  \ \mathrm{and}\
  b_k\rightarrow b
  \ \mathrm{in}\ \mathbb{R},\  \mathrm{uniformly\ on\ }[0,T],\ as\ k\rightarrow \infty.
\end{equation}
\noindent Then
\eqref{quasirelfondv^nk},
\eqref{ak-a_bk-b}, and
\eqref{eq:conv-Y_k-Y} entail
\begin{equation}\label{quasirelfondv}
  \begin{split}
    e^{-r(T-t)}& \E \left[ \langle\tilde{\lambda},
      \left(Y ^{t,y,\tilde{u}}\right)_0(T)\rangle \right] -
      v(t,y)
    \\
    &=\E
    \int_t^T  \left[ e^{-r(s-t)}\left(
      -\frac{ \left( \langle B^*\nabla  v (t,Y ^{t,y,\tilde{u}}(s))-\tilde{\beta},(1,1)\rangle \right) ^2}{4(\gamma_0+\gamma_1)}-
      \langle\tilde{\alpha},(Y ^{t,y,\tilde{u}})_0(s)\rangle\right.\right.\\
    &\left.\left.\phantom{\frac{ -\left( \langle B^*\nabla v^{k}
            (t,Y ^{t,y,\tilde{u}}(s))-\tilde{\beta},(1,1)\rangle
          \right) ^2}{4(\gamma_0+\gamma_1)}} +
      \<B\tilde {u}(s),\nabla
      v(s,Y ^{t,y,\tilde{u}}(s))\>\right)ds \right],
  \end{split}
\end{equation}
or
\begin{equation*}
  \begin{split}
    v(t,y)
    =&
    e^{-r(T-t)}
    \E \left[ \langle\tilde{\lambda},
      \left(Y ^{t,y,\tilde{u}}\right)_0(T)\rangle \right]\\
    &+\mathbb{E}
     \left[ \int_t^T
       e^{-r(s-t)}
        \left( \frac{ \left( \langle B^*\nabla  v (t,Y ^{t,y,\tilde{u}}(s))-\tilde{\beta},(1,1)\rangle \right) ^2}{4(\gamma_0+\gamma_1)}+
          \langle\tilde{\alpha},(Y ^{t,y,\tilde{u}})_0(s)\rangle\right.\right.\\
      &\left.\left.
          \phantom{
            \frac{\left(\tilde{\beta}B\right)^2}{(\gamma_0)}
          }
          -
          \<B\tilde {u}(s),\nabla
      v(s,Y ^{t,y,\tilde{u}}(s))\>
           \right) ds
         \right].
  \end{split}
\end{equation*}
Finally, adding and subtracting
\begin{equation*}
  \mathbb{E} \left[
    \int_t^T e^{-r(s-t)} \left(  \langle\tilde{\beta},\tilde{u}(s)\rangle+
      \langle\begin{bmatrix}
        \gamma_0&0\\
        0&\gamma_1
      \end{bmatrix}\tilde{u}(s),\tilde{u}(s)\rangle
    \right) ds
     \right]
\end{equation*}
we get
\begin{equation*}
  \begin{split}
    v(t,y)
    =&
    \mathcal{J}(t,y;\tilde{u})
    +\mathbb{E}
     \left[ \int_t^T
       e^{-r(s-t)}
        \left( \frac{ \left( \langle B^*\nabla  v (t,Y ^{t,y,\tilde{u}}(s))-\tilde{\beta},(1,1)\rangle \right) ^2}{4(\gamma_0+\gamma_1)}+
          \langle\tilde{\alpha},(Y ^{t,y,\tilde{u}})_0(s)\rangle\right.\right.\\
      &\left.\left.
          \phantom{
            \frac{\left(\tilde{\beta}B\right)^2}{(\gamma_0)}
          }
          -
          H_{CV}(s,B^*\nabla v(s,Y^{t,y,\tilde{u}},\tilde{u}(s))
           \right) ds
         \right].
  \end{split}
\end{equation*}

\qed

We can now pass to prove a verification theorem i.e.\
a sufficient condition of optimality given in term of the
 solution $v$ of the HJB equation.

\begin{theorem}
  \label{teorema controllo}
Let Assumption \ref{ip-concrete} hold true.
Let $v$ be in (\ref{sol-forma-special}), with $a$ and $b$ given respectively by \eqref{eq:am-sol} and \eqref{eq:bm-sol}, solution to the HJB equation
\eqref{eq:HJVcalVH_particolare}.
Then the following holds.
\begin{enumerate}
 \item For all $(t,y)\in [0,T]\times H$ we have
$v(t,y) \ge V(t,y)$, where $V$ is the value function
defined in
\eqref{optprV}.
\item
  Let $t\in [0,T],y\in H$.
  If $u^*$ is as in   \eqref{ustar}, and if $\tilde{u}^*(s)\coloneqq (u^*(B^*a(s)),u^*(B^*a(s)))$, $s\in[t,T]$,
  then
  the pair $(\tilde{u}^*,Y^{t,y,\tilde{u}^*})$ is optimal for the control problem
  \eqref{optprV},
  and $V(t,y)=v(t,y)=\mathcal{J}(t,y;\tilde{u}^*)$.
\end{enumerate}
\end{theorem}
\dim The first statement follows directly by \eqref{relfond} due to the
positivity of the integrand.
Concerning the second statement, we immediately see that, when $\tilde{u}=\tilde{u}^*$,
\eqref{relfond} becomes
$v(t,y)=\mathcal{J}(t,y;\tilde{u}^*)$.
Since we know that, for any admissible control $\tilde{u}=(u,z)\in \mathcal{U}\times \mathcal{U}$ with $u=z$,
\begin{equation*}
   \mathcal{J}(t,y;\tilde{u})\leq V(t,y) \leq v(t,x),
 \end{equation*}
the claim immediately follows.
\qed

\subsection{Equivalence with the original problem in the LQ case}
\label{SS:back}

To find the solution of the original problem in the LQ case we need to apply Proposition \ref{pr:newFG}, i.e. to prove that the optimal control
in the original LQ case is deterministic. This is the subject of next proposition.

\begin{proposition}\label{jensen}
  Condition
  \eqref{supDiag-a} is verified.
\end{proposition}
\dim
Let $t\in [0,T],y\in H$.
Let
$\tilde{u}=(u,z)\in \mathcal{U}$, with $z(s)=\mathbb{E}[u(s)]$ for $s\in [t,T]$.
Let $\tilde{u}_\mathbb{E}= (\mathbb{E}[u],z)$.
Then
\begin{equation*}
  \tilde{u}_\mathbb{E}\in
  \big\{
 \tilde{u}=(u,z)\in \mathcal{U}\times \mathcal{U},\
     \mathrm{and}\ u=z\ \mathrm{deterministic}
  \big\}.
\end{equation*}
Notice, by
\eqref{eq-mild},
that
\begin{equation}\label{YEY}
  \mathbb{E} \left[
    Y^{t,y,\tilde{u}}
  \right]
  =
  \mathbb{E} \left[
    Y^{t,y,\tilde{u}_\mathbb{E}}
  \right] .
\end{equation}
Then
\begin{equation*}
  \begin{split}
    \mathcal{J}(t,y;\tilde{u})=& \E \bigg[ \int_t^T e^{-r(s-t)} \left(
      \langle\tilde{\alpha},Y^{t,y,\tilde{u}}_0(s) \rangle
      -\langle\tilde{\beta},\tilde{u}(s)\rangle -
      \langle\begin{bmatrix}
        \gamma_0&0\\
        0&\gamma_1
      \end{bmatrix}\tilde{u}(s),\tilde{u}(s)\rangle \right)
    ds\\
    &\hspace{260pt}+\langle\tilde{\lambda},(Y^{t,y,\tilde{u}})_0(T)\rangle  \bigg]\\
    =& \E \bigg[ \int_t^T e^{-r(s-t)} \left(
      \langle\tilde{\alpha},Y^{t,y,\tilde{u}_\mathbb{E}}_0(s) \rangle
      -\langle\tilde{\beta},\tilde{u}_\mathbb{E}(s)\rangle -
      \langle\begin{bmatrix}
        \gamma_0&0\\
        0&\gamma_1
      \end{bmatrix}\tilde{u}(s),\tilde{u}(s)\rangle \right)
    ds\\
    &\hspace{260pt}+\langle\tilde{\lambda},(Y^{t,y,\tilde{u}_\mathbb{E}})_0(T)\rangle  \bigg]\\
    \leq& \textrm{(by Jensen's inequality)}\\[5pt]
    \leq & \E \bigg[ \int_t^T e^{-r(s-t)} \left(
      \langle\tilde{\alpha},Y^{t,y,\tilde{u}_\mathbb{E}}_0(s) \rangle
      -\langle\tilde{\beta},\tilde{u}_\mathbb{E}(s)\rangle -
      \langle\begin{bmatrix}
        \gamma_0&0\\
        0&\gamma_1
      \end{bmatrix}\tilde{u}_\mathbb{E}(s),\tilde{u}_\mathbb{E}(s)\rangle
    \right)
    ds\\
    &\hspace{260pt}+\langle\tilde{\lambda},(Y^{t,y,\tilde{u}_\mathbb{E}})_0(T)\rangle
    \bigg],
  \end{split}
\end{equation*}
which implies
  \eqref{supDiag-a}.

  \begin{corollary}
    Let $f,g$ be as in \eqref{lq-case}.  Let
    $t\in [0,T],x\in \mathbb{R}$.  If $u^*$ is as in \eqref{ustar},
    with $(x,x)$ in place of $y_0$, then $u^*(B^*a(s))$ is optimal for
    $\overline{V}(t,x)$.
  \end{corollary}
  \proof The statement is a straightforward consequence of
  \eqref{eq:equivVcalV}, Proposition~\ref{pr:newFG},
  Theorem~\ref{teorema controllo}.\qed


\bigskip

\textbf{Data availibility} Data sharing not applicable to this article as no datasets were generated or analysed during
the current study

\medskip

\textbf{Conflict of interest} There are no conflict of interest.

\medskip

\textbf{Ethical approval} We do not work with any empirical data. For this reason, we are not aware of any ethical issues that could arise within this article.


\begin{thebibliography}{99}

\bibitem{BFFGGEB} {R.Boucekkine, G. Fabbri, S. Federico, F. Gozzi,
\textit{A dynamic theory of spatial externalities,} Games and Economic Behavior, Elsevier, vol. 132 (C), (2022), 133-165.}

\bibitem{CarmonaDelarueBook} R. Carmona, F. Delarue, \textit{Probabilistic theory of mean field games with applications. I}
Probab. Theory Stoch. Model., 83
Springer, Cham, 2018, xxv+713 pp.

\bibitem{Cosso-et-al2022} A. Cosso, F. Gozzi, I. Kharroubi, H.  Pham, M. Rosestolato,  \textit{Optimal control of path-dependent McKean-Vlasov SDEs in infinite-dimension},
Ann. Appl. Probab. 33 (2023), no. 4, 2863-2918.

\bibitem {DP1}G. Da Prato and J. Zabczyk,\textit{\ Stochastic equations in infinite dimensions. Second Edition.}
Encyclopedia of Mathematics and its Applications 152,
Cambridge University Press, \textit{2014.}

\bibitem{deFeo} F. de Feo, \textit{Stochastic optimal control problems with delays in the state and in the control via viscosity solutions and an economical application}, arXiv:2308.14506.

\bibitem{FabbriGozziSwiech}
G. Fabbri, F. Gozzi, A. Swiech
\textit{Stochastic Optimal Control in Infinite Dimensions:
Dynamic Programming and HJB Equations}. Springer 2017.

\bibitem {GM}F. Gozzi and C. Marinelli, \textit{Stochastic optimal control of delay equations arising in advertising models}.
    Stochastic partial differential equations and applications - VII,  133-148, Lect. Notes Pure Appl. Math., 245, Chapman $\&$ Hall/CRC, Boca Raton, FL, 2006.


\bibitem{GMSJOTA} F. Gozzi, C. Marinelli, S. Savin \textit{On controlled linear diffusions with delay in a model of optimal advertising under uncertainty with memory effects}.
J. Optim. Theory Appl. 142 (2009), no. 2, 29--321.

\bibitem{FGFM-I} F. Gozzi, F. Masiero \textit{Stochastic Optimal Control with Delay in the Control, I: solving the HJB equation through partial smoothing.} SIAM J. Control Optim. 55 (2017), no. 5, pp. 2981-3012.


\bibitem{FGFM-II} F. Gozzi, F. Masiero \textit{Stochastic Optimal Control with Delay in the Control, II: Verification Theorem and Optimal Feedback Controls.} SIAM J. Control Optim. 55 (2017), no. 5, pp. 3013-3038.


\bibitem{FGFM-III} F. Gozzi, F. Masiero, \textit{Stochastic Control Problems with Unbounded Control Operators: solutions through generalized derivatives}, SIAM J. Control Optim. 61 (2023), no. 2, 586-619.

\bibitem{grosset}
L. Grosset, B. Viscolani, \textit {Advertising for a new product introduction: a stochastic
              approach},
   Top 12 (2004), no. 1, 149-167.

   \bibitem{hartl} R. F. Hartl, \textit{Optimal dynamic advertising policies for hereditary processes},
J. Optim. Theory Appl.43,
(1984),
no. 1, 51-72.

\bibitem{gw} C. Marinelli, \textit{The stochastic goodwill problem},
European J. Oper. Res. 176 (2007), no. 1, 389-404.

 \bibitem{MoPh} M. Motte and H. Pham
\textit{Optimal bidding strategies for digital advertising}, arXiv:2111.08311

\bibitem{NA} M. Nerlove, J. K. Arrow, \textit{Optimal Advertising Policy Under Dynamic Conditions}, Economica, 29 (1962), 129-142.

\bibitem{sethi-VWcomp} A. Prasad, S.P. Sethi, \textit{Competitive advertising under uncertainty: a stochastic differential game approach},
J. Optim. Theory Appl. 123 (2004), no. 1, 163-185.

\bibitem{sethi}
G. Feichtinger, R. Hartl, S. Sethi, S.,
  \textit{Dynamical {O}ptimal {C}ontrol {M}odels in {A}dvertising: {R}ecent
              {D}evelopments},
Management Sci., 40, (1994) 195-226.


\bibitem {VK}R. B. Vinter and R. H. Kwong, \textit{\ The infinite time quadratic control problem for linear systems
with state and control delays: an evolution equation approach}, SIAM J. Control Optim., 19
(1):139-153, 1981.

\bibitem{VW} M. L. Vidale, H. B. Wolfe, \textit{An operations-research study of sales response toadvertising}, Operations Res., 5, (1957) {370--381} .

\bibitem{YongZhou99} J. Yong and X.Y. Zhou
 \textit{Stochastic Control and Hamilton-Jacobi-Bellman equations.} Applications of Mathematics (New York), 43. Springer-Verlag, New York, 1999.
\end{thebibliography}
\end{document}